\documentclass[journal, a4paper]{IEEEtran}
\usepackage{amssymb}
\usepackage[utf8]{inputenc} 
\usepackage[T1]{fontenc}    
\usepackage{hyperref}       
\usepackage{url}            
\usepackage{booktabs}       
\usepackage{amsfonts}       
\usepackage{nicefrac}       
\usepackage{microtype}      
\usepackage{lipsum}		
\usepackage{graphicx}
\usepackage{natbib}
\usepackage{doi}
\usepackage{amsmath}
\usepackage{physics}
\usepackage{wrapfig}
\usepackage{bm}
\usepackage{algorithm,algpseudocode}
\usepackage{xcolor}
\usepackage{soul}
\usepackage{forest}
\usepackage{tcolorbox}
\usepackage{tikz-cd}
\usepackage{enumitem}

\begin{document}
\title{Accelerating Low-Frequency Convergence for Limited-Angle DBT via Two-Channel Fidelity in PDHG}					

\author{Taro Iyadomi,
        Ricardo Parada, Anna Kim, Lily Jiang, Emil Sidky, and William Chang 
\thanks{William Chang was with the Department of Mathematics, UCLA, CA, 90007, e-mail: chang314@g.ucla.edu.}
\thanks{Emil Sidky was with the Department of Radiology, University of Chicago, IL 60637}
\thanks{Manuscript received \today.}}

\maketitle

\begin{abstract}
Reconstruction in limited-angle digital breast tomosynthesis (DBT) suffers from slow convergence of low spatial-frequency components when using weighted data-fidelity terms within primal–dual optimization. We introduce a two-channel fidelity strategy that decomposes the sinogram residual into complementary low-pass and high-pass bands using square-root Hanning (Hann$^{1/2}$) filter families, each driven by an independent $\ell_2$-ball constraint and dual update in the PDHG (Chambolle–Pock) algorithm with He–Yuan predictor–corrector relaxation. By assigning a larger dual step-size and slightly looser tolerance to the low-frequency channel, the method delivers stronger per-iteration correction to the near-DC band without violating global PDHG stability. Experiments on a 2D digital breast phantom across multiple resolutions demonstrate that the two-channel approach yields $19\%$--$61\%$ RMSE improvement over the single-channel baseline, with larger gains at coarser discretizations where problem conditioning is more favorable, supporting more balanced spectral convergence in clinically realistic limited-angle regimes.
\end{abstract}


\section{Introduction}




Limited angular range (LAR) scanning has been widely studied in the literature, but only recently have discrete-to-discrete models demonstrated stable inversion at clinically relevant scan arcs. Building on the observation that directional sparsity improves LAR recovery, \cite{zhang2021directional} introduced a globally convex formulation using directional total variation (DTV), enabling robust optimization-based reconstruction.

We introduce a sparsity-constrained optimization framework for LAR DBT that combines directional TV penalties with pixel/voxel sparsity, enabling practical parameterization when true DTV values are unknown. To characterize the regime in which the proposed two-channel fidelity approach is most effective, we systematically vary image discretization to modulate problem conditioning and quantify the resulting reconstruction benefit.


\section{Preliminary}
We consider limited angular range (LAR) digital breast tomosynthesis (DBT) in a discrete-to-discrete setting. Let $f$ be the discretized attenuation image, $g$ is the measured sinogram data, and $X$ is the system matrix that encodes fan-beam (2D) or cone-beam (3D) projection. The idealized data model is
\begin{equation}
    g = Xf
\end{equation}
In the DBT geometry, the detector spans the in-plane directions $(x,y)$, while the $z$-direction is perpendicular to the (fixed) detector plane and is referred to as the depth direction. At clinically typical scanning arcs, the resulting linear system is generally underdetermined at useful (e.g., isotropic) discretizations, motivating the use of convex regularization.

\subsection{Directional Finite Differences and Directional Total Variation}
Let $\partial_x, \partial_y, \partial_z$ denote finite-difference operators approximating derivatives along the corresponding axes. The associated directional total-variation (DTV) seminorms are $\|\partial_x f\|_1$, $\|\partial_y f\|_1$, and $\|\partial_z f\|_1$ (in 2D studies, the $\partial_y$ term is omitted). Distinguishing in-plane and depth directional derivatives is particularly important in LAR/DBT because the data support and conditioning can differ substantially across directions.

\subsection{Constrained Sparsity-regularized Reconstruction}
\cite{zhang2021directional} showed that accurate image reconstruction from LAR data can be achieved using non-smooth convex optimization with directional total variation (DTV). Their key problem minimizes a least-squares (LSQ) data-fidelity term subject to separate $\ell_1$ bounds on the in-plane and depth finite-difference derivatives, forming DTV-constrained LSQ (DTV-LSQ):
\begin{equation}
    \begin{split}
        f^* = \arg\min_f \; \frac{1}{2}\|g - Xf\|_2^2 
\\
\text{s.t.} \quad
\|\partial_x f\|_1 \le \gamma_x, \;
\|\partial_z f\|_1 \le \gamma_z.
    \end{split}
\end{equation}

Constraining DTV is more effective than the standard total variation (TV) constraint based on gradient sparsity. DTV constraints are convenient for phantom studies where true DTV values are computable, but for real digital breast tomosynthesis (DBT) data the bounds must be estimated. The present work reformulates the problem with improved parameterization and design choices to support isotropic-resolution reconstruction in DBT.

\section{The Single Channel Approach }
In this section, we describe the single channel approach proposed by \cite{sidky2025accurate}. Let $X$ be the discrete projection operator for a limited-angle DBT scan, $g$ the measured sinogram, and $f \ge 0$ the unknown volume. The paper adopts a \emph{data-discrepancy constrained, sparsity-regularized} model:
\begin{equation}
\label{eq:model}
\begin{split}
    \min_{f \ge 0}\;\; 
\alpha_x \|\partial_x f\|_1 + \alpha_y \|\partial_y f\|_1 + \alpha_z \|\partial_z f\|_1 + \beta \|f\|_1
\\
\text{s.t.} \quad 
\big\|R[c]\,(g - X f)\big\|_2 \le \varepsilon \sqrt{\mathrm{size}(g)},
\end{split}
\end{equation}
where $R[c]$ is a square-root Hanning (Hann$^{1/2}$) filter along detector channels with cutoff parameter $c$; $\partial_\bullet$ are finite differences (directional TV). The resulting problem is convex but non-smooth (indicator constraint, $\ell_1$ penalties), and is solved with a primal-dual hybrid gradient (PDHG / Chambolle--Pock) method with a He--Yuan predictor--corrector (relaxation) using $\rho \approx 1.75$ \cite{he2012convergence}. The limited-angle regime reports running $\sim 2000$ PD iterations to reach accurate convergence.

There are two additional design choices emphasized in the paper:
(i) using $R[c]$ \emph{inside the data-constraint} to de-emphasize low spatial-frequency inconsistencies; 
(ii) including an optional voxel-sparsity term $\beta\|f\|_1$ in addition to DTV. 

\paragraph{Notation for PDHG} With $A := R[c]\,X$ and the residual $r := g - X f$, the data-constraint is $ \|A f - R[c]\,g\|_2 \le \varepsilon\sqrt{|g|}$. The PDHG dual variable $y$ lives in the sinogram space; the dual step is a projection onto an $\ell_2$ ball; the primal step applies the adjoint $A^\top = X^\top R[c]^\top$ followed by the proximal of the sum of $\ell_1$ terms and the indicator of $f\ge 0$.

\section{Why low spatial frequencies converge slowly}

The data-fidelity term can be written as $\tfrac{1}{2}\|W(g - Xf)\|_2^2$, where $W \equiv R[c]$. In the detector-frequency domain, this weighting strongly attenuates low spatial frequencies. Under a linearized analysis of the PDHG update, the correction applied to a residual Fourier mode at spatial frequency $k$ scales approximately with the squared singular value $s_k^2$ of the weighted system operator $WX$. In limited-angle DBT, the singular values $s_k$ are already small for low $k$ due to the severe ill-conditioning of $X$. The additional Hann$^{1/2}$ weighting $W$ further shrinks these low-$k$ singular values by design, creating a large spectral gap between frequency bands.

This disparity in singular values destabilizes convergence behavior. Because PDHG must satisfy the global stability condition $\tau \sigma \|WX\|^2 < 1$ (Theorem 1 of \cite{chambolle2011first}), the allowable step sizes are dictated by the \emph{largest} singular values, which reside in the mid- and high-frequency bands. Step sizes that stabilize those bands therefore provide insufficient drive for near-DC modes, yielding extremely weak contraction factors for low spatial frequencies. Increasing the global step sizes to accelerate DC correction risks divergence, since PDHG cannot apply frequency-selective step-size amplification in its standard form. Consequently, mid- and high-frequency image components converge quickly, while low spatial frequencies decay slowly and dominate the tail of the iterations. This persistent low-frequency drift reveals a fundamental mode-dependent convergence bottleneck that motivates split-frequency or multi-dual fidelity strategies.

\section{Two-channel fidelity for spectrally balanced PDHG}

Single-dual weighted fidelity in limited-angle tomographic inversion entangles spatial-frequency bands with sharply different conditioning, which can destabilize convergence when the optimizer attempts to correct near-DC error using a step size dictated by high-frequency modes. To retain the noise robustness of the Hann$^{1/2}$ weighting while providing a stronger corrective signal to low spatial frequencies, we decompose the sinogram residual into complementary spectral channels and enforce independent consistency constraints:
\begin{equation}
\label{eq:twoconstraint}
\begin{split}
    \min_{f \ge 0}\;
\alpha_x \|\partial_x f\|_1 + \alpha_y \|\partial_y f\|_1 + \alpha_z \|\partial_z f\|_1 + \beta \|f\|_1\\
\text{s.t.}\;
\begin{cases}
\|R_{\mathrm{hi}}(g - Xf)\|_2 \le \varepsilon_{\mathrm{hi}}\sqrt{|g|},\\
\|R_{\mathrm{lo}}(g - Xf)\|_2 \le \varepsilon_{\mathrm{lo}}\sqrt{|g|}.
\end{cases}
\end{split}
\end{equation}

We select both filters from the square-root Hanning family. The high-frequency channel preserves the original cutoff, $R_{\mathrm{hi}} = \text{Hann}^{1/2}(c)$, maintaining the fidelity weighting used in prior work. The low-frequency channel uses a substantially narrower cutoff, $R_{\mathrm{lo}} = \text{Hann}^{1/2}(c_{\mathrm{lo}})$ with $c_{\mathrm{lo}} \ll c$, forming a complementary pair that partitions spectral correction responsibilities while keeping $R_{\mathrm{hi}}^2 + R_{\mathrm{lo}}^2$ close to identity over the inversion-relevant passband. The low-frequency constraint tolerance is chosen slightly looser than the high-frequency channel to allow mild model mismatch at DC while still driving large-scale bias reduction.

\begin{algorithm}[t]
\caption{PDHG with split fidelity channels}
\begin{algorithmic}[1]
\State \textbf{Inputs:} $g, X, R_{\mathrm{hi}}, R_{\mathrm{lo}}, \alpha_{x,y,z}, \beta, \varepsilon_{\mathrm{hi}}, \varepsilon_{\mathrm{lo}}, \sigma_{\mathrm{hi}}, \sigma_{\mathrm{lo}}, \tau, \rho$
\State Initialize $f^0$, $\bar f^0 = f^0$, $y_{\mathrm{hi}}^0 = y_{\mathrm{lo}}^0 = 0$
\For{$k = 0,1,2,\dots$}
    \State $y_{\mathrm{hi}}^{k+1} \gets \operatorname{Proj}_{B_{\mathrm{hi}}}\!\big(y_{\mathrm{hi}}^k + \sigma_{\mathrm{hi}} R_{\mathrm{hi}}(g - X\bar f^k)\big)$
    \State $y_{\mathrm{lo}}^{k+1} \gets \operatorname{Proj}_{B_{\mathrm{lo}}}\!\big(y_{\mathrm{lo}}^k + \sigma_{\mathrm{lo}} R_{\mathrm{lo}}(g - X\bar f^k)\big)$
    \State $f^{k+1} \gets \operatorname{prox}_{\tau\big(\alpha_x\|\partial_x \cdot\|_1 + \alpha_y\|\partial_y \cdot\|_1 + \alpha_z\|\partial_z \cdot\|_1 + \beta\|\cdot\|_1\big) + \mathbb{I}\{f \ge 0\}}$ $\Big(f^k - \tau X^\top(R_{\mathrm{hi}}^\top y_{\mathrm{hi}}^{k+1} + R_{\mathrm{lo}}^\top y_{\mathrm{lo}}^{k+1})\Big)$
    \State $\bar f^{k+1} \gets f^{k+1} + \rho(f^{k+1} - f^k)$ \Comment{He–Yuan relaxation}
\EndFor
\end{algorithmic}
\end{algorithm}

\paragraph{Spectral rationale}
In primal–dual gradient schemes, the convergence rate of each Fourier mode is governed by the corresponding singular values of the forward operator. In limited-angle systems, high-frequency modes are well represented and produce comparatively large singular values, while low spatial frequencies induce much smaller gains. The original Hann weighting amplifies this gap by further suppressing near-DC singular values. When a single dual is used, the step size must be set to stabilize the high-frequency spectrum, which unintentionally results in an under-driven update for the low-frequency band. The optimizer therefore applies only a weak correction to DC error each iteration, but a strong correction to high frequencies, creating oscillatory or stalled behavior in the low band.

By assigning each channel an independent dual step size, the method reduces cross-band interference. The low-frequency dual ascent can safely be made larger, providing a stronger contraction factor for small $k$, without compromising the stability required to converge the high-frequency channel. The proximal update then combines both adjoint-filtered corrections while preserving non-negativity. This split-fidelity design yields a more spectrally uniform correction profile and mitigates frequency-induced convergence bottlenecks in limited-angle tomographic inversion.

\section{Experimental Methodology}

To evaluate the proposed two-channel fidelity approach, we conduct inverse problem studies using a 2D LAR scanning configuration matching that of Sidky et al.~\cite{sidky2025accurate}. The geometry uses 25 fan-beam projections over a 50$^\circ$ arc. The source is positioned 50~cm from the isocenter and 100~cm from the detector, which consists of 1024 bins spanning a length that inscribes a 10~cm $\times$ 10~cm field of view.

\subsection{Digital Phantom}

The phantom used for these studies is a 2D digital breast slice following the methodology of~\cite{sidky2025accurate}. Glandular tissue regions are created by applying a threshold to power-law noise, with small bright specks added to simulate microcalcifications~\cite{reiser2010task}. The test image is constructed as a weighted combination of three tissue components: adipose ($0.5\times$), fibroglandular ($1.0\times$), and calcification ($2.0\times$) layers.

To characterize how the two-channel fidelity benefit depends on problem conditioning, we test three array sizes: $512 \times 512$ pixels (0.2~mm pixel width), $256 \times 256$ pixels (0.4~mm pixel width), and $128 \times 128$ pixels (0.8~mm pixel width). Because the LAR inverse problem has a fixed number of measurements (25 views $\times$ 1024 detector bins $= 25{,}600$ ray integrals), coarser discretizations yield fewer unknowns and thus better-conditioned systems. This experimental design allows us to isolate the relationship between conditioning and the efficacy of the split-frequency fidelity strategy. Following~\cite{sidky2025accurate}, projection data are generated using the same pixel grid as the reconstruction, which removes discretization error and focuses the comparison on algorithmic differences. All studies use noiseless projection data to isolate the convergence behavior of the optimization algorithms from noise-related effects.

\subsection{Algorithm Configuration}

Both reconstruction methods---single-channel and two-channel---are configured with matched regularization parameters to ensure fair comparison. Following the resolution-dependent parameterization of~\cite{sidky2025accurate}, the DTV weight parameter is set to $\alpha = 1.7$ for $512 \times 512$, $\alpha = 1.9$ for $256 \times 256$, and $\alpha = 1.95$ for $128 \times 128$. The pixel sparsity parameter is set to $\beta = 5.0$ for the $512 \times 512$ case and $\beta = 10.0$ for the $256 \times 256$ and $128 \times 128$ cases.

The PDHG algorithm is run for 500 iterations using the He--Yuan predictor-corrector scheme~\cite{he2012convergence} with relaxation parameter $\rho = 1.75$. For the single-channel baseline, the data-discrepancy constraint employs a square-root Hanning filter with cutoff parameter $c = 4.0$. The two-channel method uses complementary filters: a low-pass channel with cutoff $c_{\mathrm{lo}} = 8.0$ and a high-pass channel with cutoff $c_{\mathrm{hi}} = 4.0$. The low-frequency channel employs an amplified dual step size with $\sigma_{\mathrm{lo}} = 4.0 \times \sigma_{\mathrm{hi}}$ and a relaxed tolerance $\varepsilon_{\mathrm{lo}} = 1.25 \times \varepsilon_{\mathrm{hi}}$, enabling stronger per-iteration correction to near-DC components while maintaining global PDHG stability~\cite{chambolle2011first}. Step sizes are determined via power iteration to satisfy the convergence condition $\tau \sigma \|WX\|^2 < 1$. Our implementation is publicly available at \url{https://github.com/taroii/ct}.

\section{Results}

We now present quantitative and visual comparisons of the single-channel and two-channel methods. Table~\ref{tab:rmse} summarizes reconstruction error across three image resolutions.

\begin{table}[h]
\renewcommand{\arraystretch}{1.3}
\caption{Image RMSE Comparison Across Resolutions}
\label{tab:rmse}
\centering
\begin{tabular}{lccc}
\toprule
Resolution & Single-Channel & Two-Channel & Improvement \\
\midrule
$512 \times 512$ & 0.1625 & 0.1317 & 19.0\% \\
$256 \times 256$ & 0.0984 & 0.0767 & 22.1\% \\
$128 \times 128$ & 0.0332 & 0.0128 & 61.4\% \\
\bottomrule
\end{tabular}
\end{table}

The two-channel method achieves consistently lower reconstruction error across all tested resolutions. Notably, the improvement magnitude increases substantially with coarser discretization: from 19.0\% at $512 \times 512$ to 61.4\% at $128 \times 128$. This trend directly reflects the relationship between problem conditioning and the efficacy of split-frequency fidelity. At coarser resolutions, the ratio of measurements to unknowns improves ($25{,}600$ measurements versus $16{,}384$ unknowns at $128 \times 128$, compared to $262{,}144$ unknowns at $512 \times 512$), yielding a better-conditioned inverse problem. Under improved conditioning, the low-frequency error components---which the two-channel approach specifically targets---become the dominant source of residual error, and the amplified dual step size for the low-frequency channel can provide proportionally greater correction. Conversely, at finer resolutions where severe ill-conditioning introduces errors across all frequency bands, the targeted low-frequency correction constitutes a smaller fraction of the total improvement opportunity.

\begin{figure}[h]
    \centering
    \includegraphics[width=1\linewidth]
    {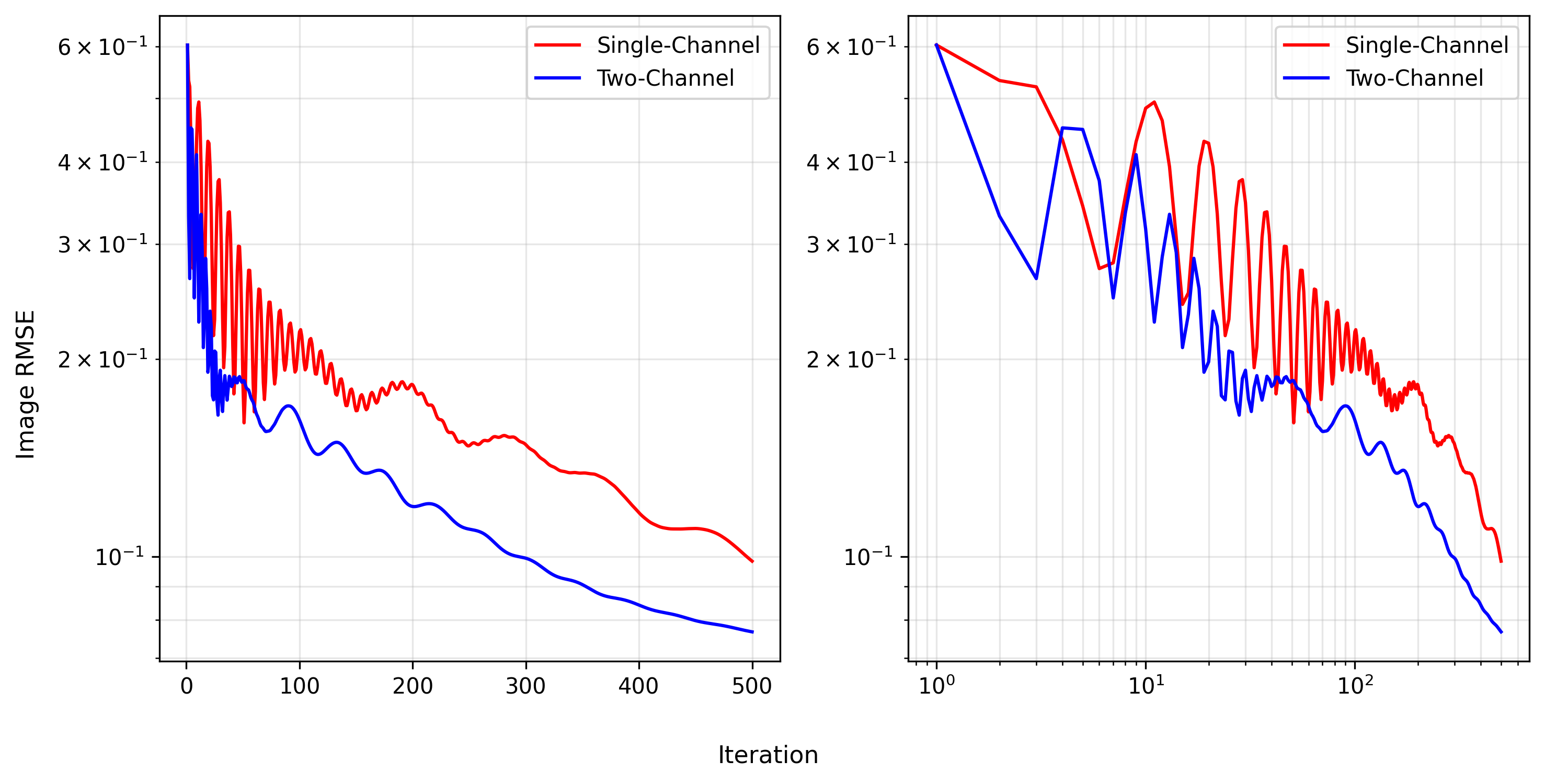}
    \caption{Image RMSE convergence for $256 \times 256$ reconstruction. Left: linear iteration scale. Right: log-log scale. The two-channel method (blue) achieves lower error and reduced early-iteration oscillation compared to the single-channel baseline (red).}
    \label{fig:convergence}
\end{figure}

Figure~\ref{fig:convergence} shows the convergence behavior for the $256 \times 256$ case. Both methods exhibit decrease in image RMSE, but the two-channel approach achieves lower error throughout the iteration process. The single-channel method displays pronounced oscillations in the early iterations that persist until approximately iteration 200, whereas the two-channel method stabilizes more rapidly after an initial transient. This difference is consistent with the theoretical motivation: by assigning an independent dual step size to the low-frequency channel, the two-channel approach reduces cross-band interference that would otherwise cause the optimizer to oscillate as it attempts to balance corrections across frequency bands with disparate conditioning.

\begin{figure}[h]
    \centering
    \includegraphics[width=0.8\linewidth]{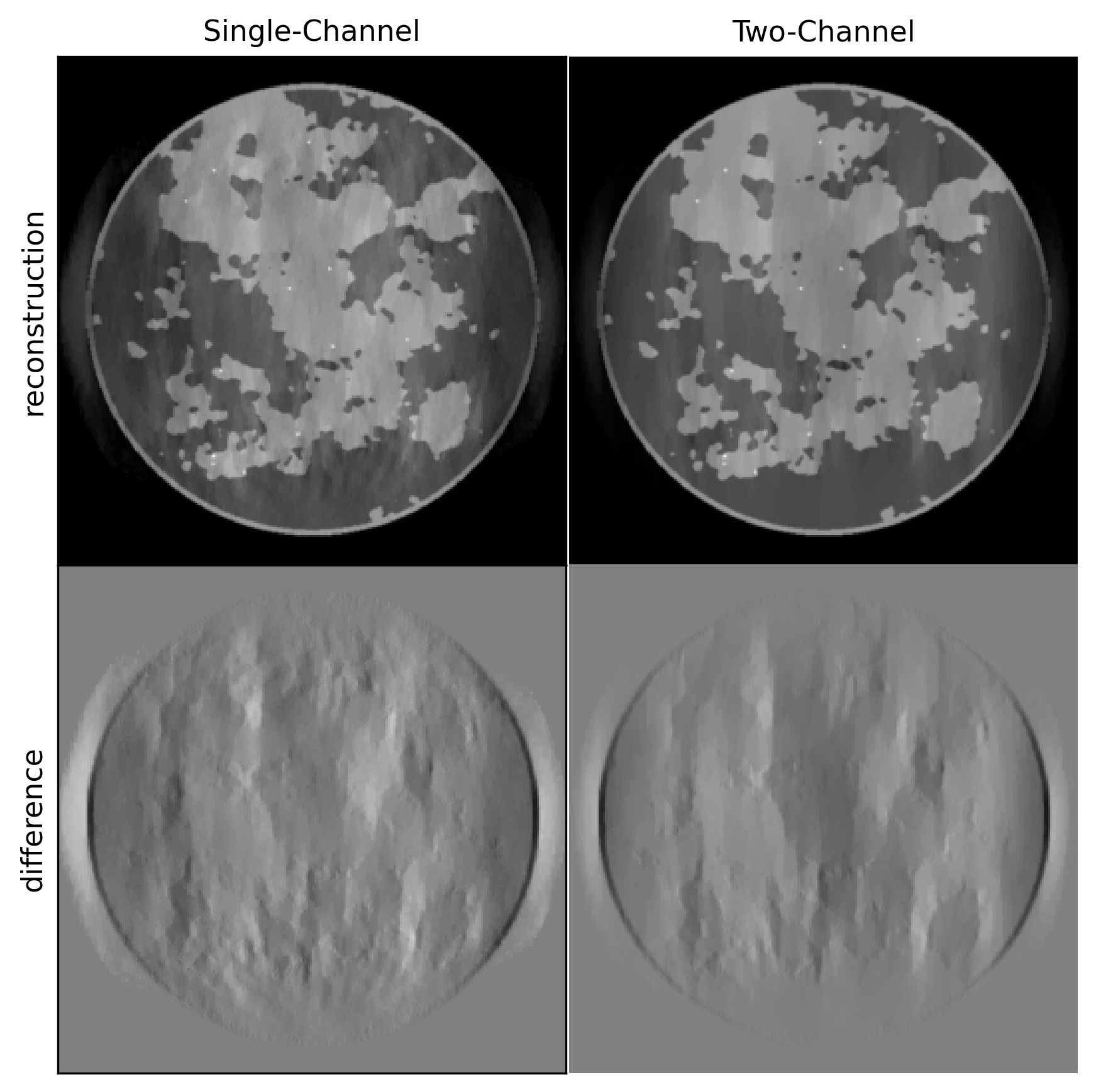}
    \caption{Reconstruction comparison for $256 \times 256$ images. Top row: single-channel (left) and two-channel (right) reconstructions. Bottom row: corresponding difference images (reconstruction minus ground truth), displayed with grayscale range $[-0.75, 0.75]$. The two-channel result shows reduced large-scale intensity artifacts.}
    \label{fig:recon_comparison}
\end{figure}

Figure~\ref{fig:recon_comparison} presents the reconstructed images and difference maps for the $256 \times 256$ case at 500 iterations. Both methods recover the gross phantom structure, including the glandular tissue distribution and microcalcification specks. The difference images reveal that the single-channel method retains more residual error distributed throughout the reconstruction, particularly in the form of large-scale intensity variations characteristic of unconverged low-frequency components. The two-channel reconstruction exhibits more uniform error reduction, with the residual errors appearing more localized to fine structural boundaries.

These results are consistent with the spectral analysis of Section 4: frequency-separated data fidelity enables more balanced convergence across the spatial frequency spectrum. The benefit is most pronounced when the underlying inverse problem is better conditioned, as occurs at coarser discretizations where the number of unknowns more closely approaches the number of measurements.

\section{Conclusion}

We analyzed the fundamental causes of slow and unstable convergence in limited-angle digital breast tomosynthesis reconstruction when using a single, globally weighted fidelity term in primal–dual optimization. The square-root Hanning weighting creates a large singular-value gap between high and near-DC modes. Because PDHG step sizes must stabilize the largest gains (high frequencies), the low-frequency channel is under-driven, leading to stalled or unstable convergence at DC.

To address this, we proposed a two-channel fidelity design that routes low and high spatial-frequency residuals through complementary square-root Hanning filters, each paired with an independent dual variable and step size. This strategy preserves noise robustness while reducing cross-band step-size interference, providing stronger corrective feedback to the ill-conditioned low-frequency channel without violating global PDHG stability. The resulting formulation better equalizes spectral corrections, mitigates low-band oscillations, and directly targets the most ill-conditioned modes in limited-angle inversion.

Experiments across multiple image resolutions reveal that the two-channel benefit scales with problem conditioning: improvements range from 19\% at $512 \times 512$ (severely underdetermined) to 61\% at $128 \times 128$ (moderately underdetermined). This conditioning dependence suggests the method is most impactful when low-frequency errors dominate the residual, as occurs when the inverse problem is sufficiently well-posed for high-frequency components to converge reliably.

These findings support the broader conclusion that accounting for spectral conditioning gaps in reconstruction optimizers for limited-angle tomography can help avoid unstable convergence tails. As future work, we believe that multi-dual or frequency-split fidelity designs can offer a principled and practical path to improving stability and convergence uniformity in other types of tomographic reconstruction.

\bibliography{references}
\bibliographystyle{abbrvnat}

\end{document}